\documentclass[12pt]{article}

\usepackage{amsmath}
\usepackage{amsfonts}
\usepackage{amssymb}

\newtheorem{proposition}{Proposition}[section]
\newtheorem{theorem}{Theorem}[section]
\newtheorem{lemma}[theorem]{Lemma}
\newtheorem{corollary}[theorem]{Corollary}
\newtheorem{remark}[theorem]{Remark}

\newtheorem{definition}{Definition}

\def\phi{{\varphi}}

\DeclareSymbolFont{AMSb}{U}{msb}{m}{n}
\DeclareMathSymbol{\N}{\mathbin}{AMSb}{"4E}
\DeclareMathSymbol{\Z}{\mathbin}{AMSb}{"5A}
\DeclareMathSymbol{\R}{\mathbin}{AMSb}{"52}
\DeclareMathSymbol{\Q}{\mathbin}{AMSb}{"51}
\DeclareMathSymbol{\I}{\mathbin}{AMSb}{"49}
\DeclareMathSymbol{\C}{\mathbin}{AMSb}{"43}

\begin{document}
\title{Convergence of the alternating split Bregman algorithm in infinite-dimensional Hilbert spaces}

\author{{Amir Moradifam\footnote{Department of Mathematics, University of Toronto, Toronto, Ontario, Canada M5S 2E4. E-mail: amir@math.toronto.edu. The author is supported by a MITACS Postdoctoral Fellowship. }
\qquad Adrian Nachman\footnote{Department of Mathematics and the
Edward S. Rogers Sr. Department of Electrical and Computer
Engineering, University of Toronto, Toronto, Ontario, Canada M5S
2E4. E-mail: nachman@math.toronto.edu. The author is supported in part by an NSERC Discovery Grant.}}}
\date{\today}

\smallbreak \maketitle

\begin{abstract}

We prove results on weak convergence for the alternating split Bregman algorithm in infinite dimensional Hilbert spaces. We also show convergence of an approximate split Bregman algorithm, where errors are allowed at each step of the computation. To be able to treat the infinite dimensional case, our proofs focus mostly on the dual problem. We rely on Svaiter's theorem on weak convergence of the Douglas-Rachford splitting algorithm and on the relation between the alternating split Bregman and Douglas-Rachford splitting algorithms discovered by Setzer. Our motivation for this study is to provide a convergent algorithm for weighted least gradient problems arising in the hybrid method of imaging electric conductivity from interior knowledge (obtainable by MRI) of the magnitude of one current. 

\end{abstract}

\section{Introduction}
Split Bregman and alternating split Bregman algorithms were proposed by Goldstein and Osher \cite{GO} for solving problems of the form
\[\min_{u} f(d)+g(u) \ \ \hbox{subject to} \ \ d=Du,\]
where $D\in \R^{m\times n}$ is a linear transform acting from $\R^n$ to $\R^m$ and $f, g$ are  convex functions. These algorithms have been shown to be very successful in various PDE based image restoration approaches and compressed sensing \cite{GO,ZBO}. In particular, the alternating split Bregman algorithm is very efficient for large scale $l^1$-norm minimization and TV minimization problems problems \cite{GO,ZBO}. The convergence of the split Bregman algorithm was proved in \cite{GO}. Later in \cite{COS} and \cite{Setzer, Setzer1}  the authors independently proved convergence of the alternating split Bregman algorithm in finite dimensional Hilbert spaces. In this paper we shall prove weak convergence results for the alternating split Bregman algorithm in infinite dimensional Hilbert spaces. Our motivation for this study is to provide a convergent algorithm for weighted least gradient problems arising in the hybrid method of imaging electric conductivity from interior knowledge (obtainable by MRI) of the magnitude of one current. Our proof relies on a recent result of Svaiter \cite{S} about weak convergence of the Douglas-Rachford splitting algorithm.

Let $H_1$ and $H_2$ be real Hilbert spaces and consider the minimization problem
\begin{equation*}
(P) \hspace{1cm}\min_{u \in H_1} \{ g(u)+f(L u)\},
\end{equation*}
where $L: H_1\rightarrow H_2$ is a bounded linear operator and both functions $g:H_1\rightarrow \R \cup \{\infty\}$ and $f:H_2\rightarrow \R \cup \{\infty\}$ are proper, convex and lower semi-continuous.  The problem $(P)$ can be written as a constrained minimization problem
\begin{equation}
\min_{u\in H_1, d\in H_2}  g(u)+f(d) \ \ \ \ \hbox{subject to} \ \ \ \ Lu=d,
\end{equation}
which leads to an unconstrained problem:
\begin{equation}
\min_{u\in H_1, d\in H_2}  g(u)+f(d) +\frac{\lambda}{2} \|Lu-d\|^2. 
\end{equation}
To solve the above problem, Goldstein and Osher \cite{GO} introduced the split Bregman method: 
\begin{eqnarray}\label{SB}
(u^{k+1},d^{k+1})&=&\hbox{argmin}_{u\in H_1, d\in H_2} \{g(u)+f(d)+\frac{\lambda}{2}\parallel b^k+Lu -d\parallel^2_{2}\}, \\
b^{k+1}&=& b^k+Lu^{k+1}-d^{k+1}.\nonumber
\end{eqnarray}
Yin et al \cite{Yinetal} (see also \cite{TW}) showed that the split Bregman algorithm can be viewed as an augmented Lagrangian algorithm \cite{H, P, R}.  Since the joint minimization problem (\ref{SB}) in both $u$ and $d$ could sometimes be hard to solve exactly, Goldstein and Osher \cite{GO} proposed the following algorithm for solving the problem $(P)$. \\

\textbf{Alternating split Bregman algorithm:} \\

Initialize  $b^0$ and  $d^0$. 
For  $k\geq1$:\\
\begin{enumerate}
\item Find a minimizer $u^k$ of 
\begin{equation}
I^{k}_1(u):= g(u)+\frac{\lambda}{2}\parallel b^{k-1}+Lu-d^{k-1}\parallel^{2}_2,
\end{equation}
on $H_1$. 
\item  Find the minimizer $d^k$ of
\begin{equation}
I_2^k(d):=f(d)+\frac{\lambda}{2}\parallel b^{k-1}+Lu^{k}-d\parallel^{2}_2,
\end{equation}
on $H_2$. 
\item Let $ b^{k}=b^{k-1}+Lu^{k}-d^{k}$. 
\end{enumerate}

Cai, Osher, and Shen \cite{COS} proved that if the primal problem $(P)$ has a unique solution then the sequence $u^k$ in the above algorithm will converge to the minimizer of $(P)$.  Independently, Setzer \cite{Setzer} showed that the alternating split Bregman algorithm coincides with the Douglas-Rachford splitting algorithm applied to the dual problem 

\begin{equation*}
(D) \hspace{.5cm}-\min_{b \in H_2} \{ g^*(-L^*b)+f^*(b)\},
\end{equation*}
and proved the convergence of the alternating split Bregman algorithm in finite-dimensional Hilbert spaces (See \cite{Setzer1}, Proposition 1 and Theorem 4).

 We note that for general L, the functional $I^k_{1}(u)$ in the above algorithm may not have a minimizer in $H_1$ and therefore the alternating split Bregman algorithm  may not be well defined. We thus make the following definition, which will be used throughout this paper. 

\begin{definition}
We say that the alternating split Bregman algorithm is well defined if  for all $k\geq 1$ the functionals $I^k_{1}(u)$ and $I^k_{2}(d)$ have minimizers in $H_1$ and $H_2$, respectively.  
\end{definition}

Note that $I^k_{2}(u)$ is strictly convex and coercive in $H_2$. Hence if $f$ is weakly lower semi continuous, then  $I^k_{2}(d)$ will have a unique minimizer in $H_1$. The following proposition  provides a sufficient condition for  existence of a minimizer of $I^k_{1}(u)$. 
\begin{proposition} \label{prop3^*}
Let $H_1$ and $H_2$ be two Hilbert spaces, $g:H_1\rightarrow \R \cup \{\infty\}$ be a proper, convex, and lower semi-continuous function. Assume $L: H_1\rightarrow H_2$ is a bounded linear operator. If $L^*L:H_1\rightarrow H_1$ is surjective, then for every $c\in H_2$ the functional 
\begin{equation}\label{I}
I(u)= g(u)+\frac{\lambda}{2}\parallel Lu+c\parallel^{2}_2,
\end{equation} 
has a unique minimizer on $H_1.$ 
\end{proposition}

We include a proof of Proposition 1.1 in Section 2 of this paper. We are now ready to state our main theorems. Below is a special case of a more general result (Theorem 2.6) that we will prove in Section 2. 
We let $L^*$ denote the Hermitian adjoint of L. Throughout the paper, we make the usual identification of a dual of a Hilbert space $H$ with $H$ itself. 

\begin{theorem}\label{main-special}
Let $H_1$ and $H_2$ be two Hilbert spaces (possibly infinite dimensional) and assume that both primal $(P)$ and dual $(D)$ problems have optimal solutions and that $L^*L:H_1\rightarrow H_1$ is surjective. Then the alternating split Bregman algorithm is well defined and the sequences  $\{b^k\}_{k\in N}$  and $\{d^k\}_{k\in N}$ converge weakly to some  $\hat{b}$ and $\hat{d}$, respectively. Moreover    $\lambda \hat{b}$ is a solution of the dual problem, $\{d^k-L u^{k+1}\}_{k\in N}$ converges strongly to zero, and
\begin{equation}\label{conv-est}
\sum_{k=0}^{\infty}||d^k-L u^{k+1}||^2_{2}<\infty.
\end{equation}
Furthermore, there exists a unique $\hat{u}\in H_1$ such that $L\hat{u}= \hat{d}$ and $\hat{u}$ is a solution of the primal problem $(P)$. In particular $\{u^k\}_{k\in N}$  has at most one weak cluster point $\bar{u}=\hat{u}$. 
\end{theorem}

Under additional assumptions on the functionals $f$ and $g$ we will also prove the following result. 

\begin{theorem}\label{COS}
Let $H_1$ and $H_2$ be two Hilbert spaces (possibly infinite dimensional). Suppose that alternating split Bregman is well defined,  $f$ is continuous, and both $f$ and $g$ are weakly lower semi continuous. Then the sequences  $\{b^k\}_{k\in N}$  and $\{d^k\}_{k\in N}$ converge weakly to some  $\hat{b}$ and $\hat{d}$, respectively. Moreover    $\lambda \hat{b}$ is a solution of the dual problem $(D)$, $\{d^k-L u^{k+1}\}_{k\in N}$ converges strongly to zero, and
\begin{equation}\label{conv-estCOS}
\sum_{k=0}^{\infty}||d^k-L u^{k+1}||^2_{2}<\infty.
\end{equation}
Furthermore
\begin{equation}\label{energy}
\lim_{k\rightarrow \infty} f(Lu^k)+g(u^k)=\min_{u \in H_1} f(Lu)+g(u),
\end{equation}

every weak cluster point of  $\{u^k\}_{k\in N}$  is a solution of the problem $(P)$, and  $L(\hat{u})=\hat{d}$ for every cluster point $\hat{u}$ of $\{u^k\}_{k\in N}$. In particular if (P) has a unique solution then the sequence $\{u^k\}_{k\in N}$  has at most one weak cluster point which is a solution of the problem (P). 

\end{theorem}

Comparing Theorem \ref{COS} with Theorem \ref{main-special} (and with our more general result Theorem \ref{main}) one can see the effect of the operator $L$ on the convergence behaviour  of the alternating split Bregman algorithm.  Indeed if $L$ is injective then Theorem \ref{main}  guarantees that the sequence $\{u^k\}_{k\in N}$ has at most one weak cluster point without assuming uniqueness of minimizers of the problem (P), continuity of $f$, or weak lower semi continuity of the functionals $f, g$.  When $L$ is injective and the problem $(P)$ has more than one solution then, depending on the initial values of $b^0$ and $d^0$, the alternating split Bregman algorithm may converge to different solutions of the primal problem.  On the other hand when $L$ is not injective and the primal problem (P) has a unique solution then Theorem \ref{COS}  guarantees that $\{u^k\}_{k\in N}$ has at most one weak cluster point while  Theorem \ref{main} only says $L^{-1}(\hat{d})$ contains a solution of the primal problem, which is a weaker conclusion.

In many applications $L^*L$ is surjective and Theorem \ref{main-special} guarantees the convergence of the alternating split Bregman algorithm. For instance in weighted least gradient problems $H_1=H^1_0(\Omega)$, $H_2=(L^2(\Omega))^n$, and $L=\nabla u$. It is easy to check that $L^*:(L^2(\Omega))^n \rightarrow H^{-1}_0(\Omega)$ is surjective. Recently, our group studied the problem of recovering an isotropic conductivity
from the interior measurement of the magnitude of one current density field
\cite{MNT, NTT07, NTT08}. We showed that the conductivity is uniquely determined by the magnitude of the current generated by imposing a given boundary voltage. Moreover the corresponding voltage potential is the unique minimizer of the infinite-dimensional minimization problem 
\begin{equation}\label{min_prob}
u=\hbox{argmin} \{\int_{\Omega}|J||\nabla v|: v \in H^{1}(\Omega),
\ \ v|_{\partial \Omega}=f\},
\end{equation}
where $|J|$ is the magnitude of the current density vector field generated by imposing the voltage $f$ on the boundary of the connected bounded region $\Omega \subset R^n$, $n\geq 2$. 

The results presented in this paper lead to a convergent split Bregman algorithm for computing the unique minimizer of the least gradient problem (\ref{min_prob}). The details will be presented in a forthcoming paper \cite{MNTim}, along with a number of successful numerical experiments for recovering the electric conductivity.

\section{Convergence of the alternating split Bregman algorithm }

Recall that, by Fenchel duality \cite{Rb}, the dual problem corresponding to the problem $(P)$ can be written as 
\begin{equation*}
(D) \hspace{.5cm}-\min_{b \in H_2} \{ g^*(-L^*b)+f^*(b)\}.
\end{equation*}
Let $v(P)$ and $v(D)$ be the optimal values of the primal (P) and dual problem (D), respectively. Weak duality always holds, that is $v(P)\geq v(D)$ \cite{FPA, Rb}. To guarantee strong duality, i.e., the equality $v(P)=v(D)$ together with existence of a solution to the dual problem, several regularity conditions are available in the literature (see \cite{FPA}, Chapter 7).  In particular, if
\[\exists x \in \hbox{dom}(g) \cap \hbox{dom} (foL) \ \ \hbox{such that $f$ is continuous at $Lx$}, \]
or
\[0\in \hbox{int} (L\hbox{dom}(g)-\hbox{dom}(f)),\]
then strong duality holds. In this paper  we will always assume that both primal and dual problems have optimal solutions, and at least one of the above conditions is satisfied. Then by Rockafellar-Fenchel duality \cite{Rb}, if $b$ is any solution of the dual problem, then the \textit{entire }set of solutions of the primal problem is obtained as
 \begin{equation}\label{RF}
 \partial g^*(-L^* b) \cap L^{-1}\partial f^*(b). 
 \end{equation}
The above representation of solutions of the primal problem is the key for our understanding of the alternating split Bregman algorithm. To explain this, we first prove the following simple lemma.

\begin{lemma}\label{lem}
Let $g:H_1\rightarrow \R \cup \{\infty\}$ and assume $L: H_1\rightarrow H_2$ is a bounded linear operator.  Assume one of the following conditions hold:
\begin{enumerate}
\item There exists a point $L\bar{u}$ where $g^*$ is continuous and finite. 
\item $L^*:H_2\rightarrow H_1$ is surjective.  
\end{enumerate}
Then
\begin{equation}\label{eqval0}
\partial (g^*o(-L^*))(b) =-L\partial g^*(-L^* b),
\end{equation}
for all $b\in H_2$. 
 In particular if  $L^{-1}(d) \neq \emptyset$ for some $d\in H_2$, then
\begin{equation}\label{eqval}
 -d \in \partial (g^*o(-L^*))(b) \ \  \Leftrightarrow \ \ L^{-1}d \cap  \partial g^*(-L^* b)\neq \emptyset.
\end{equation}
\end{lemma}
{\bf Proof.} If 1) holds, then (\ref{eqval0}) follows from Proposition 5.7 in \cite{ET}.  Now assume $L^*$ is surjective and let $-L(u)=-d \in \partial (g^*o(-L^*))(b) $.  Then 
\begin{eqnarray*}
g^*(-L^*(c))-g^*(-L^*(b)) \geq \langle -Lu, c-b\rangle= \langle u, -L^*(c-b)\rangle,
\end{eqnarray*}
for all $c \in H_2$. Since $L^*$ is surjective, $u\in \partial g^*(-L^* b)$ and consequently $L^{-1}d \subseteq  \partial g^*(-L^* b)$. 

Now assume  $u\in \partial(-L^*(b))$. Then 
\begin{eqnarray*}
g^*(-L^*(c))-g^*(-L^*(b)) \geq  \langle u, -L^*(c-b)\rangle=\langle -Lu, c-b\rangle=\langle-d, c-b\rangle,
\end{eqnarray*}
for all $c\in H_2$. Hence $ -d \in \partial (g^*o(-L^*))(b)$. \hfill $\Box$\\

By above lemma,  if one can find $\hat{d}\in L(H_1)$ such that 
\[\hat{d}\in \partial f^*(\hat{b}) \ \ \hbox{and} \ \ -\hat{d} \in \partial (g^*o(-L^*))(\hat{b}),\]
for some solution $\hat{b}$ of the dual problem, then 
\begin{equation}\label{heart}
L^{-1}(\hat{d}) \subseteq \partial g^*(-L^* b) \cap L^{-1}\partial f^*(b).
\end{equation} 
Consequently, by Rockafellar-Fenchel duality,  every $\hat{u}\in L^{-1}\hat{d}$ will be a solution of the primal problem $(P)$.  One can find $\hat{d} \in H_2$ satisfying the above conditions by solving the inclusion problem 
\begin{equation}\label{doualPartial}
 0 \in \partial (g^*o(-L^*))(\hat{b})+\partial f^*(\hat{b}). 
 \end{equation} 
 
Indeed if $\hat{b}\in H_2$ is a solution of the dual problem and $\hat{d}\in \partial f^*(\hat{b})$, then $-\hat{d} \in \partial (g^*o(-L^*))(\hat{b})$ and (\ref{heart}) follows from Lemma \ref{lem}.  This is summarized in the following lemma. 

\begin{lemma}\label{heartLemma} Let $g:H_1\rightarrow \R \cup \{\infty\}$ and $f:H_2\rightarrow \R \cup \{\infty\}$ be two proper lower semi-continuous convex functions and assume that hypothesis of Lemma \ref{lem} hold. Also suppose both primal (P) and dual (D) problems have optimal solutions and let $\hat{b}$ be an arbitrary  solution of the dual problem,  then every $\hat{u} \in L^{-1}(\partial f^*(\hat{b}))$ is a solution of the primal problem (P). 

\end{lemma}
 
We thus focus on computing a solution of the problem (\ref{doualPartial}). This can be written in the form of an inclusion problem
 \begin{equation}\label{inclusion}
 0\in A(\hat{b})+B(\hat{b}),
 \end{equation}
 where $A:=\partial (g^*o(-L^*))$ and $B:=\partial f^*$ are maximal monotone operators on $H_2$. If both primal and dual problems have optimal solutions, then the above inclusion problem has at least one solution. Therefore, if we can find a solution $\hat{b}$  of the problem (\ref{inclusion}) as well as $\hat{d}\in B(\hat{b})$, then by Lemma \ref{heartLemma} every $\hat{u}\in L^{-1}(\hat{d})$ will be a solution of the primal problem (D).  
The Douglas-Rachford splitting  method in Convex Analysis provides precisely such a pair $(\hat{b},\hat{d})$. Following this route leads to the alternating split Bregman algorithm. Indeed it is shown by Setzer in \cite{Setzer} that the alternating split Bregman algorithm for the primal problem (P) coincides with Douglas-Rachford spliting algorithm applied to (\ref{inclusion}) (see Theorem \ref{setzer}).

We now explain this in more detail.Let $H$ a real Hilbert space and $A, B: H\rightarrow 2^{H}$ be two maximal monotone operators. For a set valued function $P: H\rightarrow 2^{H}$, let $J_P$ to be its resolvent i.e.,
\[J_{P}=(Id+P)^{-1}.\]
It is well known that sub-gradient of convex, proper, lower semi-continuous functions are maximal monotone \cite{Rb} and if $P$ is maximal monotone then $J_P$ is single valued.

Lions and Mercier \cite{LM}  showed that for any general maximal monotone operators $A,B$ and any initial element $x_0$ the sequence
Defined by the Douglas-Rachford recursion: \begin{equation}\label{iter:DR}
x_{k+1}=(J_A(2J_B-Id)+Id-J_B)x_k,
\end{equation}
converges weakly to some point $\hat{x} \in \emph{H}$ such that $\hat{p}=J_B (\hat{x})$ solves the inclusion problem (\ref{inclusion}). Much more recently, Svaiter \cite{S} proved that the sequence $p_k=J_{ B}(x_k)$ also converges weakly to $\hat{p}$. 

\begin{theorem}{(\bf Svaiter \cite{S})}\label{DRS}
Let $H$ be a Hilbert space and $A, B: \emph{H}\rightarrow 2^{\emph{H}}$ be maximal monotone operators and assume that a solution of (\ref{inclusion}) exists. Then, for any initial elements $x_0$ and $p_0$ and any $\lambda>0$, the sequences $p_{k}$ and $x_k$ generated by the following algorithm
\begin{eqnarray}
x_{k+1}&=&J_{\lambda A}(2p_k-x_k)+x_k-p_k \nonumber \\
p_{k+1}&=&J_{\lambda B}(x_{k+1}),
\end{eqnarray}
converges weakly to some $\hat{x}$ and $\hat{p}$ respectively. Furthermore, $\hat{p}=J_{\lambda B}(\hat{x})$  and $\hat{p}$ satisfies
\[0\in A(\hat{p})+B(\hat{p}).\]
\end{theorem}

To apply Douglas-Rachford splitting algorithm one needs to evaluate the resolvents $J_{\lambda A}(2p_k-x_k)$ and $J_{\lambda B}(x_{k+1})$ at each iteration.  To evaluate the resolvents we are led to find  minimizers of $I^k_1(u)$ and $I^k_2(d)$ in the alternating split Bregman algorithm. Indeed if we let 
\[A=\partial (g^*o(-L^*)),  \ \ B=\partial f^*, \ \ x^0=\lambda(b^0+d^0), \ \ \hbox{and} \ \ p^0=\lambda b^0.\]  
Then the resolvents $J_{\lambda A}(2p_k-x_k)$ and $ J_{\lambda B}(x_{k+1})$ can be computed as follows 
\[J_{\lambda A}(2p_k-x_k)=\lambda(b^k+Lu^{k+1}-d^k),\]
and 
\[J_{\lambda B}(x_{k+1})=\lambda (b^k+Lu^{k+1}-d^{k+1}),\]
where $u^{k+1}$ and $d^{k+1}$ are minimizers of $I^k_1(u)$ and $I^k_2(d)$, respectively (see \cite{Setzer, Setzer1} for a proof). The following theorem gives the precise relation between the sequences generated by the alternating split Bregman algorithm and those generated by the Douglas-Rachford splitting algorithm.

\begin{theorem}{(\bf Setzer \cite{Setzer})}\label{setzer}
The Alternating Split Bregman Algorithm coincides with Douglas-Rachford splitting  algorithm applied to (D) with $A:=\partial (g^*o(-L^*))$ and $B:=\partial f^*$, where 
\begin{equation}
x_k=\lambda(b^k+d^k), \ \ \ \ p_k=\lambda b^k, \ \ k\geq 0. \\
\end{equation}
\end{theorem}
\vspace*{0.5cm}

The operator $T:=J_A(2J_B-Id)+Id-J_B$ is known to be firmly non-expansive, i.e.,  $T=\frac{1}{2}Id+\frac{1}{2}R$, with $R$ satisfying:
\[\parallel Rx -Ry \parallel \leq \parallel x-y\parallel \ \ \ \ \hbox{for all} \ \ x,y \in H.\]

We will need the following lemma in our convergence proof. 
\begin{lemma}\label{nonexpansice-lemma}  If $T:H\rightarrow H$ is a firmly non-expansive operator and $x_{k+1}=T(x_k)$ with $x_0 \in H$, then
\[ \parallel x_{k+1}-\hat{x}\parallel ^2+\parallel x_{k+1}-x_k\parallel ^2 \leq \parallel x_{k}-\hat {x} \parallel ^2.\]
\end{lemma}
{\bf Proof. } Since $T$ is firmly non-expansive, $R=2T-Id$ is a non-expansive operator.

 Hence
\[\parallel R x_k- R \hat{x}\parallel ^2 =-\parallel x_k- \hat{x} \parallel ^2+2\parallel T x_k -T\hat{x}\parallel^2-2 \parallel (Id-T)x_k -(Id-T)\hat{x}\parallel^2. \]
Therefore we have
\[\frac{1}{2}(\parallel x_k- \hat{x} \parallel ^2 -\parallel R x_k- R \hat{x}\parallel ^2 )=\parallel x_k- \hat{x} \parallel ^2-\parallel x_{k+1} -\hat{x}\parallel^2- \parallel x_{k+1}-x_k\parallel^2.\]
Since $R$ is non-expansive, the left hand side of the above inequality is non-negative, and this completes the proof. \hfill $\Box$\\

Now we are ready to prove our main theorem. 

\begin{theorem}\label{main}
Let $H_1$ and $H_2$ be two Hilbert spaces (possibly infinite dimensional) and assume that both the primal $(P)$ and the dual $(D)$ problems have optimal solutions. Suppose that (\ref{eqval0}) holds and  the alternating split Bregman algorithm is well defined.  Let $\{u^k\}_{k\in N}$, $\{d^k\}_{k\in N}$, and $\{b^k\}_{k\in N}$ be the three sequences generated by the alternating Split Bregman algorithm. Then  $\{d^k\}_{k\in N}$, and $\{b^k\}_{k\in N}$  converge weakly to some  $\hat{d}$, and $\hat{b}$, respectively.  Moreover $\lambda \hat{b}$ is a  solution of the dual problem, the sequence  $\{d^k-L u^{k+1}\}_{k\in N}$ converges strongly to zero, and
\begin{equation}\label{conv-est}
\sum_{k=0}^{\infty}||d^k-L u^{k+1}||^2_{2}<\infty.
\end{equation}
Furthermore, $L^{-1}(\hat{d})$ contains a solution $\hat{u}$ of the primal problem $(P)$.  In particular if $L$ is injective, then  $\{u^k\}_{k\in N}$  has at most one weak cluster point $\bar{u}=\hat{u}$. 

\end{theorem}

{\bf Proof.} The weak convergence of the sequences  $d^k$, and $b^k$  follows from Theorems \ref{DRS} and \ref{setzer}. To prove the estimate (\ref{conv-est}), let $T=J_A(2J_B-Id)+Id-J_B$. Since $T$ is firmly non-expansive, by Lemma \ref{nonexpansice-lemma} we have
\begin{eqnarray}
\parallel x_{k+1}-\hat{x}\parallel ^2+\parallel x_{k+1}-x_k\parallel ^2 \leq \parallel x_{k}-\hat {x} \parallel ^2,
\end{eqnarray}
where $\hat{x}$ is the weak limit of $x_k$ with $T(\hat{x})=\hat{x}$.  By the above inequality, we have
\begin{equation}
\sum_{k=0}^{\infty}\parallel x_{k+1}-x_{k}\parallel^2 <\infty. 
\end{equation}
Now observe that
\[x_k-x_{k-1}=\lambda ( b^{k+1}+d^{k+1}-b^k-d^{k})=\lambda (Lu^{k+1}-d^k),\]
and hence (\ref{conv-est}) follows.  

By Theorem \ref{setzer}  and Theorem \ref{DRS},  $\hat{p}=\lambda \hat{b}$ is a minimizer of the dual problem and $J_{\lambda \partial f^*}  ( \lambda (\hat{d}+\hat{b}))=\lambda \hat{b}$.  Therefore
\[\lambda \hat{b}+\lambda \partial f^* (\lambda \hat{b})=\lambda (\hat{d}+\hat{b}) \ \ \Leftrightarrow \ \ \hat{d} \in \partial f^*(\lambda \hat{b}) \ \ \Leftrightarrow -\hat{d} \in \partial (g^*o(-L^*))(\lambda \hat{b}).\]

By Lemma \ref{lem} there exists $\hat{u} \in H_1$ such that  $\hat{u}\in \partial g^*(-L^*(\hat{p}))$ and $L(\hat{u})=\hat{d}$.  Therefore 
\[\hat{u}\in \partial g^*(-L^*(\hat{p})) \cap L^{-1}(\partial f^*(\hat{p})).\]
Since $\hat{p}$ is a minimizer of the dual problem, it follows from the Fenchel-Rockafellar duality theorem that $\hat{u}$ is a minimizer of the primal problem.  If $L$ is injective then  the sequence $u^k$ has at most one weak cluster point $\bar{u}$ and necessarily $\bar{u}=\hat{u}$. The proof is now complete. \hfill $\Box$\\

{\bf Proof of Theorem \ref{COS}:} Since $f$ is continuous by Theorem 4.1 in \cite{ET} the dual problem (D) has an optimal solution. Thus it follows from an argument similar to that of Theorem \ref{main}  that the sequences $b^k$ and $d^k$ weakly converge to some $\hat{b}$ and $\hat{d}$ where $\lambda \hat{b}$ is a solution of the dual problem and (\ref{conv-estCOS}) holds.  In particular if $\hat{u}$ is a weak cluster point of $u^k$ then $L(\hat{u})=\hat{d}$. To prove (\ref{energy}), we can now use the argument of Cai, Osher, and Shen in the proof of Theorem 3.2 in \cite{COS} Let $\hat{u}$ be a solution of the problem $(P)$ and set $\hat{d}=L\hat{u}$,  $\hat{p} \in \partial f(\hat{d})$, and $\hat{b}=\frac{\hat{p}}{\lambda}$.  Define 
\[u^k_e=u^k-\hat{u}, \ \ d^k_e=d^k-\hat{d}, \ \ b^k_e=b^k-\hat{b}.\]
Then, as in \cite{COS},

\begin{eqnarray}\label{COSequ}
&&\frac{\lambda}{2}(\parallel b^0_e\parallel^2-\parallel b^{K+1}_e\parallel^2 +\parallel d^0_e\parallel^2-\parallel b^{K+1}_e\parallel^2)\nonumber \\
&=&\sum_{k=0}^{K}\langle \partial g(u^{k+1})-\partial g(\hat{u}), u^{k+1}-\hat{u}\rangle + \sum_{k=0}^{K}\langle \partial f (d^{k+1})-\partial f(\hat{d}), d^{k+1}-\hat{d}\rangle \nonumber \\
&+&\frac{\lambda}{2} \left( \sum_{k=0}^{K} \parallel Lu^{k+1}_e-d^{k+1}_e\parallel^2 +\sum_{k=0}^{K} \parallel Lu^{k+1}_e-d^{k}_e\parallel^2\right).
\end{eqnarray}
Now (\ref{energy}) follows from  (\ref{COSequ}) as in the proof of  Theorem 3.2 in \cite{COS}.  Finally since both $f$ and $g$ are weakly lower semi continuous, in view of (\ref{energy}), every weak cluster point of  $\{u^k\}_{k\in N}$ is a solution of the primal (P). \hfill $\Box$ \\

\begin{remark}
Notice that Theorem \ref{DRS} is crucial for the proof of the convergence of the sequences $b^k$ and $d^k$.\\
\end{remark}

{\bf Proof of Proposition \ref{prop3^*}:}  First note that $\hat{u}$ is a minimizer of (\ref{I}) if and only if 
\begin{equation}\label{Ip}
0\in (1/\lambda)\partial  g(\hat{u})+L^*L\hat{u}+L^*c \\ \Leftrightarrow \\  \hat{u} \in \left( (1/\lambda)\partial  g +L^*L\right)^{-1}(-L^*c).
\end{equation}
Therefore to guarantee existence of a solution of (\ref{I}) it is enough to prove that $(1/\lambda)\partial  g +L^*L$ is surjective.  Let $A:=L^*L$. Since $A$ is the subgradient of the convex lower semi-continuous functional $\parallel Lu\parallel^2$, it is a maximal monotone operator.  We claim that $A$ is $3^*-$monotone, i.e. 
\begin{equation}\label{3*}
\forall (x,x^*)\in H_1\times H_1 \ \ \ \  \sup_{(y,y^*) \in Graph (A)} \langle x-y,y^*-x^*\rangle <\infty.
\end{equation}
Since $A$ is surjective and $(y,y^*) \in Graph (A)$, 
\[y^*=A(y), \ \ \hbox{and}\ \ x^*=A(x_0) \ \ \hbox{for some} \ \ x_0 \in H_1.\]
Hence 
\begin{eqnarray*}
\langle x-y,y^*-x^*\rangle&=&\langle L(x-y),L(y-x_0)\rangle\\
&=&-\| L(x_0-y)\|^2+\langle L(x-x_0),L(y-x_0)\rangle\\
&\leq &-\| L(x_0-y)\|^2+\| L(x-x_0) \| \|L(y-x_0)\|\\
&\leq& \frac{\| L(x-x_0) \|}{4}< \infty. 
\end{eqnarray*}
Thus (\ref{3*}) holds and $A$ is $3^*-$monotone. It follows from Theorem 8 in \cite{Ma} and Lemma 2.2 in \cite{B1} that the operator  $(1/\lambda)\partial  g +L^*L$ is surjective and therefore $(\ref{I})$ has a solution.

Now notice that since $L^*L$ is surjective, $L^*$ is surjective and hence $L$ is injective. Consequently the functional $I(u)$ is strictly convex and  has a unique minimizer.  The proof is now complete. \hfill $\Box$

\begin{remark}
If $A=L^*L$ is surjective then one can show that $A$ is also invertible. To see this note that if $A$ is surjective then $L^*$ is surjective and hence $L$ is injective. Now assume $L^*L(u)=0$. Then
\[0=\langle L^*L(u), u\rangle=\|Lu \|^2.\]
Therefore $u=0$. 
\end{remark}

Theorem \ref{main-special} now follows from Theorem \ref{main}.

%\begin{remark}
%Intuitively the estimate (\ref{conv-est}) asserts that $||d^k-L u^{k+1}||^2_{2}$ decays faster than $1/k$.  One can compare this estimate with the estimate in the Split Bregman method i.e.,
%\[||d^k-L u^{k}||^2\leq \frac{C}{k},\]
%\end{remark}
%(see Theorem 2.1 in \cite{GO}). 

\section{Approximate alternating split Bregman algorithm}
In this section we show that the alternating split Bregman algorithm is stable with respect to possible errors at each step in the calculation of minimizers of $I^k_1(u)$ and $I^k_2(d)$.  The proof relies on the following theorem about the Douglas-Rachford splitting algorithm.  

\begin{theorem} \label{theo:SP}{\bf (Svaiter \cite{S})} Let $\lambda>0$, and let $\{\alpha_k\}_{k \in N}$ and $\{\beta_k\}_{k \in N}$ be sequences in a Hilbert space $H$. Suppose $0 \in\hbox{ran}(A+B)$, and $\sum_{k \in N} (\parallel \alpha _k\parallel +\parallel \beta_k\parallel)<\infty$.  Take $x_0 \in H$ and set
\begin{equation}
x_{k+1}=x_k+J_{\gamma A}(2(J_{\lambda B}x_k+\beta_n)-x_k)+\alpha_k -(J_{\lambda B}x_k+\beta_k) , \ \ k\geq 1. 
\end{equation}
Then $x_k$ and $p_k=J_{\lambda B}x_k$ converge weakly to $\hat{x} \in H$ and  $\hat{p} \in H$, respectively and $\hat{p}=J_{\lambda B}\hat{x} \in (A+B)^{-1}(0)$. \\
\end{theorem}

The proof of the above theorem in infinite-dimensional Hilbert spaces is due to Svaiter \cite{S} (see also \cite{C}). \\

\textbf{Approximate alternating split Bregman algorithm:} \\

Initialize  $b^0$ and  $d^0$. For $k\geq1$:\\
\begin{enumerate}
\item Find $u^k$ such that
\[\parallel Lu^k-Lu^k_{ex}\parallel_{H_2} \leq \alpha_k, \]
where $u^{k}_{ex} $ is a minimizer of 
\[I^{k}_1(u)= \{g(u)+\frac{\lambda}{2}\parallel b^{k-1}+Lu-d^{k-1}\parallel^{2}_2\},\]
on $H_1$. 
\item  Find $d^k$ such that 
\[\parallel d^k -d^k_{ex}  \parallel_{H_2}\leq  \beta_{k},\]
where $d^{k}_{ex}$ is the minimizer of 
\[I_2^k(d)=\hbox{argmin}_{d \in H_2} \{f(d)+\frac{\lambda}{2}\parallel b^{k-1}+Lu^{k}-d\parallel^{2}_2\},\]
on $H_2$. 
\item Let $ b^{k}=b^{k-1}+Lu^{k}-d^{k}$. 
\end{enumerate}

 By Theorem \ref{theo:SP} and an argument similar to that of Theorem \ref{main}  we can prove the following theorem about convergence of the sequences $u^k$, $d^k$, and $b^k$ produced by the above algorithm. 
 
 \begin{theorem}\label{mainPurt}
Let $H_1$ and $H_2$ be two Hilbert spaces (possibly infinite dimensional) and assume that both primal (P) and dual $(D)$ problems have optimal solutions. Suppose that (\ref{eqval0}) holds, the perturbed alternating split Bregman algorithm is well defined, and
\[\sum_{k=1} ^{\infty}(\alpha _k +\beta_k)<\infty.\]
 Let $\{u^k\}_{k\in N}$, $\{d^k\}_{k\in N}$, and $\{b^k\}_{k\in N}$ be the three sequences generated by the approximate alternating Split Bregman algorithm. Then  $\{d^k\}_{k\in N}$, and $\{b^k\}_{k\in N}$  converge weakly to some  $\hat{d}$, and $\hat{b}$, respectively.  Moreover $\lambda \hat{b}$ is a  solution of the dual problem, the sequence  $\{d^k-L u^{k+1}\}_{k\in N}$ converges strongly to zero, and
\begin{equation*}\label{conv-est}
\sum_{k=0}^{\infty}||d^k-L u^{k+1}||^2_{2}<\infty.
\end{equation*}
Furthermore, $L^{-1}(\hat{d})$ contains a solution $\hat{u}$ of the primal problem $(P)$.  In particular if $L$ is injective, then  $\{u^k\}_{k\in N}$  has at most one weak cluster point $\bar{u}=\hat{u}$. 
\end{theorem}

If $L^*L$ is surjective, we have the following stronger result.

\begin{corollary}
Let $H_1$ and $H_2$ be two Hilbert spaces (possibly infinite dimensional) and assume that both the primal (P) and the dual $(D)$ problems have optimal solutions and that $L^*L:H_1\rightarrow H_1$ is surjective. Suppose that 
\[\sum_{k=1}^{\infty} (\alpha _k +\beta_k)<\infty.\]
Let $\{u^k\}_{k\in N}$, $\{d^k\}_{k\in N}$, and $\{b^k\}_{k\in N}$ be the three sequences generated by the approximate alternating Split Bregman algorithm, then the sequences $\{d^k\}_{k\in N}$, and $\{b^k\}_{k\in N}$  converge weakly to some  $\hat{d}$, and $\hat{b}$, respectively.  Moreover $\lambda \hat{b}$ is a  solution of the dual problem, the sequence  $\{d^k-L u^{k+1}\}_{k\in N}$ converges strongly to zero, and
\begin{equation*}\label{conv-est}
\sum_{k=0}^{\infty}||d^k-L u^{k+1}||^2_{2}<\infty.
\end{equation*}
Furthermore there exists a unique $\hat{u}\in H_1$ such that $L\hat{u}= \hat{d}$ and $\hat{u}$ is a solution of the primal problem $(P)$.  In particular, $\{u^k\}_{k\in N}$  has at most one weak cluster point $\bar{u}=\hat{u}$. 
\end{corollary}

\end{document}